\documentclass[12pt,reqno, a4paper, final]{amsart}

\usepackage{a4}                          
\usepackage[english]{babel}               
\usepackage[latin1]{inputenc}            
\usepackage{amsmath}                     
\usepackage{amssymb}                     
\usepackage{amsthm}
\usepackage{graphicx}
\usepackage{hyperref}
\usepackage{mathrsfs}
\usepackage[all]{xy}
\usepackage[T1]{fontenc}
\usepackage{url}
\usepackage{array}
\usepackage{verbatim,epsfig}

\theoremstyle{plain}

\newtheorem{thm}{Theorem}[section] 
\newtheorem{cor}[thm]{Corollary}
\newtheorem{lem}[thm]{Lemma}
\newtheorem{prop}[thm]{Proposition}
\newtheorem{defi}[thm]{Definition}

\newtheorem{obs}[thm]{Observation}

\newcommand{\GG}{{\mathbb G}}

\newcommand{\NN}{{\mathbb N}}
\newcommand{\ZZ}{{\mathbb Z}}

\newcommand{\RR}{{\mathbb R}}
\newcommand{\CC}{{\mathbb C}}

\newcommand{\TT}{{\mathbb T}}
\newcommand{\MM}{{\mathbb M}}

\newcommand{\ip}[2]{\langle {#1} \hspace{0.03cm} | \hspace{0.03cm} {#2} \rangle}

\renewcommand{\L}{{\mathscr L}}

\newcommand{\M}{{\mathscr M}}
\newcommand{\U}{{\mathscr U}}

\newcommand{\A}{{\mathscr A}}

\newcommand{\B}{{\mathscr B}}

\newcommand{\I}{{\mathbf{I}}}

\newcommand{\varps}{{\varepsilon}}

\newcommand{\htens}{\bar{\otimes}}

\newcommand{\tens}{\otimes}
\renewcommand{\Im}{{\operatorname{Im}}}
\renewcommand{\Re}{{\operatorname{Re}}}

\newcommand{\To}{\longrightarrow}

\newcommand{\red}{{\operatorname{red}}}
\newcommand{\Tor}{\operatorname{Tor}}
\newcommand{\op}{{\operatorname{{op}}}}

\newcommand{\Hom}{\operatorname{Hom}}

\renewcommand{\H}{\mathscr{H}}
\newcommand{\id}{\operatorname{id}}

\newcommand{\bet}{\beta^{(2)}}

\newcommand{\hh}{H^{(2)}}
\newcommand{\gl}{\mathfrak{gl}}
\newcommand{\g}{\mathfrak g}
\newcommand{\Alg}{{\operatorname{Alg}}}
\newcommand{\alg}{{\operatorname{alg}}}

\renewcommand{\d}{\operatorname{d \hspace{-0.065cm}}}
\newcommand{\mat}{{\mathbb M}}

\newcommand{\utens}[1]{\underset{#1}{\tens}}
\newcommand{\algutens}[1]{\underset{#1}{\odot}}
\renewcommand{\ln}{\log}

\begin{document}

\title{\texorpdfstring{$L^2$-homology for compact quantum groups}{L2-homology for compact quantum groups}}
\author{David Kyed }
\date{}
\address{David Kyed \\ Institute for Mathematical Sciences\\
  Uni\-versi\-tets\-parken 5 \\ DK-2100 Copenhagen \O \\ Denmark}

\email{kyed@math.ku.dk}
\urladdr{\url{www.math.ku.dk/~kyed}}

\begin{abstract}
A notion of $L^2$-homology for compact quantum groups is introduced, generalizing the classical notion for countable, discrete groups. If the compact quantum group in question has tracial Haar state, it is possible to define its $L^2$-Betti numbers and Novikov-Shubin invariants/capacities. It is proved that these $L^2$-Betti numbers vanish for the Gelfand dual of a compact Lie group and that the zeroth Novikov-Shubin invariant equals the dimension of the underlying Lie group. Finally, we relate our approach  to the approach of A.~Connes  and D.~Shlyakhtenko by proving that the $L^2$-Betti numbers of a compact quantum group, with tracial Haar state, are equal to the Connes-Shlyakhtenko $L^2$-Betti numbers of its Hopf $*$-algebra of matrix coefficients.
\end{abstract}

\maketitle

\section{Introduction and definitions}
The notion of $L^2$-invariants was introduced by  M.~F.~Atiyah in \cite{atiyah} in the setting of a Riemannian manifold endowed with a free, proper and cocompact action of a discrete, countable group. Later this notion was generalized by J.~Cheeger and M.~Gromov in \cite{CG} and by W.~Lück in \cite{Luck98}. The latter of these generalizations makes it possible to define $L^2$-homology and $L^2$-Betti numbers of an arbitrary topological space equipped with an arbitrary action of a discrete, countable group $\Gamma$. In particular, the $L^2$-homology and $L^2$-Betti numbers of $\Gamma$, which are defined in \cite{Luck98} using the action of $\Gamma$ on $E\Gamma$, make sense and can be expressed in the language of homological algebra as
\[
H_n^{(2)}(\Gamma)=\Tor_n^{\CC\Gamma}(\L(\Gamma),\CC) \quad \textrm{ and } \quad \bet_n(\Gamma)=\dim_{\L(\Gamma)}H_n^{(2)}(\Gamma),
\]
where $\dim_{\L(\Gamma)}(\cdot)$ is W.~Lück's generalized Murray-von Neumann dimension introduced in \cite{Luck98}. A detailed introduction to this dimension theory can  be found in \cite{Luck02}.\\ 

	Consider now a compact quantum group $\GG=(A,\Delta)$ in the sense of S.~L.~Woro\-no\-wicz; i.e.~$A$ is a unital $C^*$-algebra and $\Delta\colon A\to A\tens A$ is a unital, coassociative $*$-homomorphism satisfying a certain non-degeneracy condition. We shall not elaborate further on the notion of compact quantum groups, but refer the reader to the survey articles \cite{vandaele} and \cite{tuset} for more details.  Denote by $h$ the Haar state on $\GG$ and by $(A_0,\Delta,S,\varps)$ the Hopf $*$-algebra of matrix coefficients arising from irreducible corepresentations of $\GG$. We recall (\cite[Prop.~7.8]{vandaele}) that $h$ is faithful on $A_0$. Consider the GNS-representation $\pi_h$ of $A$ on $L^2(A,h)$ and denote by $\M$ the enveloping von Neumann algebra $\pi_h(A)''$. We then make the following definition:
\begin{defi}\label{def1}
The $n$-th $L^2$-homology of the compact quantum group $\GG=(A,\Delta)$ is defined as
\[
H_n^{(2)}(\GG)=\Tor_n^{A_0}(\M,\CC).
\]
Here $\CC$ is considered a left $A_0$-module via the counit $ \varps$ and $\M$ is considered a right $A_0$-module via the natural inclusion $\pi_h\colon A_0 \to \M$. The groups $H_n^{(2)}(\GG)$ have a natural left $\M$-module structure and when the Haar state $h$ is tracial we may therefore define the $n$-th $L^2$-Betti number of $\GG$ as
\[
\bet_n(\GG)= \dim_{\M}H_n^{(2)}(\GG),
\]
where $\dim_\M(\cdot)$ is W.~Lück's extended dimension function arising from the extension to $\M$ of the tracial Haar state $h$.
\end{defi}    
Similar to the algebraic extension of the notion of Murray-von Neumann dimension, the classical notion of Novikov-Shubin invariants was transported to an algebraic setting by W.~Lück (\cite{Luck97}) using finitely presented modules, and generalized to arbitrary modules by W.~Lück, H.~Reich and T.~Schick in \cite{LRS}. This generalization was worked out using capacities which are essentially inverses of Novikov-Shubin invariants (cf.~\cite[Section 2]{LRS}).  In particular, they defined the $n$-th capacity of a discrete countable group $\Gamma$ as $c_n(\Gamma)= c(H_n^{(2)}(\Gamma))$, the right-hand side being the capacity of the $n$-th $L^2$-homology of $\Gamma$ considered as a left module over the group von Neumann algebra $\L(\Gamma)$. Following this approach we make the following definition:
\begin{defi}\label{def2}
If $h$ is tracial we define the $n$-th capacity of $\GG$ as $c_n(\GG)=c(H_n^{(2)}(\GG))$.
\end{defi}
To justify Definition \ref{def1} and Definition \ref{def2} we prove the following:
\begin{prop}\label{liggruppe}
Let $\Gamma$ be a countable discrete group and consider the compact quantum group $\GG=(C^*_\red(\Gamma),\Delta_\red)$ where $\Delta_\red$ is defined by $\Delta_\red\lambda_\gamma=\lambda_\gamma\tens\lambda_\gamma$ and $\lambda$ denotes the left regular representation of $\Gamma$. Then $H_n^{(2)}(\GG)=H_n^{(2)}(\Gamma)$ and in particular
\[
\bet_n(\GG)=\bet_n(\Gamma) \quad \textrm{ and } \quad c_n(\GG)=c_n(\Gamma),
\]
for all $n\in \NN_0$.
\end{prop}
Here, and in what follows, $\NN_0$ denotes the set of non-negative integers.
\begin{proof}
Since the Haar state on $\GG$ is the  trace state $\tau(x)=\ip{x\delta_e}{\delta_e}$, the GNS-action of $C^*_\red(\Gamma)$ on $L^2(C^*_\red(\Gamma),\tau)$ is naturally identified with the standard action of $C^*_\red(\Gamma)$ on $\ell^2(\Gamma)$. Note that each $\lambda_\gamma$ is a one-dimensional
(hence irreducible) corepresentation of $\GG$ and that these span a dense subspace in $L^2(C^*_\red(\Gamma),h)\simeq \ell^2(\Gamma)$. It now follows from the quantum Peter-Weyl Theorem (\cite[Thm.~3.2.3]{tuset}) that the Hopf $*$-algebra of matrix coefficients coincides with $\lambda(\CC\Gamma)$ and from this we see that the counit coincides with the trivial representation of $\Gamma$. Thus
\[
H_n^{(2)}(\GG)=\Tor_n^{(C_\red^*(\Gamma))_0}(\L(\Gamma),\CC)=\Tor_n^{\CC\Gamma}(\L(\Gamma),\CC)=H_n^{(2)}(\Gamma).
\]
\end{proof}
In the following sections we shall focus on computations of $L^2$-invariants for some concrete compact quantum groups. More precisely, the paper is organized as follows: In Section \ref{zerosection} we focus on the zeroth $L^2$-Betti number and capacity and prove that if the compact quantum group in question is the Gelfand dual $C(G)$ of a compact Lie group $G$ with $\dim(G)\geq 1$, then the zeroth $L^2$-Betti number vanishes and the zeroth capacity equals the inverse of $\dim(G)$. Section \ref{comsection} is devoted to proving that also the higher $L^2$-Betti numbers of the abelian quantum group $C(G)$ vanish in the case when $G$ is a compact, connected Lie group. In Section \ref{cs-section} we prove that the $L^2$-Betti numbers of a compact quantum group  $\GG=(A,\Delta)$ with tracial Haar state $h$ are equal to the Connes-Shlyakhtenko $L^2$-Betti numbers (see \cite{CS}) of the tracial $*$-algebra $(A_0,h)$.
\vspace{0.3cm}
\paragraph{\emph{Acknowledgements.}} I wish to thank my supervisor Ryszard Nest for suggesting that I study $L^2$-invariants in the context of quantum groups, and for the many discussions and ideas about the subject along the way. Moreover, I thank the people at the Mathematics department in Göttingen for their hospitality during the early summer of 2006 where parts of the work were carried out.
\vspace{0.3cm}
\paragraph{\emph{Notation.}}
All tensor products between $C^*$-algebras occurring in the following are assumed to be minimal/spatial. These will be denoted $\tens$ while tensor products in the category of Hilbert spaces and the category of von Neumann algebras will be denoted $\htens$. Algebraic tensor products will be denoted $\odot$.

\section{\texorpdfstring{The zeroth $L^2$-invariants}{The zeroth L2-invariants}}\label{zerosection}
In this section we focus on the zeroth $L^2$-Betti number and capacity. The first aim is to prove that the zeroth  $L^2$-Betti number of a compact quantum group, whose enveloping  von Neumann algebra is a finite factor, vanishes. After that we compute the zeroth $L^2$-Betti number and capacity for Gelfand duals of compact Lie groups and finally we study the $L^2$-invariants of finite dimensional quantum groups.
\subsection{The factor case}
In this subsection we investigate the case when the enveloping von Neumann algebra is a finite factor. First a small observation.

\begin{obs}\label{two-sidet}
Let $\M$ be a von Neumann algebra and let $A_0$ be a strongly dense $*$-subalgebra of $\M$. Let  $J_0$ be a two-sided ideal in $A_0$ and denote by $J$ the left ideal in $\M$ generated by $J_0$. Then the strong operator closure $\bar{J}$ is a two-sided ideal in $\M$. Clearly $\bar{J}$ is a left ideal and because $A_0$ is dense in $\M$  we get that $xm\in \bar{J}$ whenever $x\in J$ and $m\in \M$. From this it follows that $\bar{J}$ is also a right ideal in $\M$. 
\end{obs}

The following proposition should be compared to \cite[Cor.~2.8]{CS}.
\begin{prop}\label{faktor-nul}
Let $\GG=(A,\Delta)$ be a compact quantum group with tracial Haar state $h$.  Denote by $\pi_h$ the GNS-representation of $A$ on $L^2(A,h)$ and assume that $\M=\pi_h(A_0)''$ is a factor. If $A\neq \CC$ then $\bet_0(\GG)=0$.
\end{prop}
\begin{proof}
First note that
\[
H_0^{(2)}(\GG)=\Tor_0^{A_0}(\M,\CC)\simeq \M\algutens{A_0}\CC \simeq \M/J, 
\]
where $J$  is the left ideal in $\M$ generated by $\pi_h(\ker(\varps))$. Since the counit $\varps\colon A_0\to \CC$ is a $*$-homomorphism its kernel is a two-sided ideal in $A_0$, and by Observation \ref{two-sidet} we conclude that the strong closure $\bar{J}$ is a two-sided ideal in $\M$. Since $A\neq \CC$ the kernel of $\varps$ is non-trivial and hence $\bar{J}$ is nontrivial. Any finite factor is simple (\cite[Cor.~6.8.4]{KR2}) and therefore $\bar{J}=\M$. Now note that
\[
J\subseteq \bar{J} \subseteq \overline{J}^{\alg}=\hspace{-0.5cm}\bigcap_{\underset{J\subseteq \ker(f)}{f\in \Hom_\M(\M,\M)}}\hspace{-0.5cm}\ker(f).
\]
Since $\M$ is finitely (singly) generated as an $\M$-module,  \cite[Thm.~0.6]{Luck98} implies that $\dim_\M(J)=\dim_\M(\overline{J}^\alg)$ and thus 
\[
\dim_\M(J)=\dim_\M(\bar{J})=\dim_\M(\M)=1. 
\]
Additivity of the dimension function (\cite[Thm.~0.6]{Luck98}) now yields the desired conclusion.

\end{proof}
Denote by $A_o(n)$ the free orthogonal quantum group. The underlying  $C^*$-algebra $A$ is the universal, unital $C^*$-algebra generated by $n^2$ elements $\{u_{ij}\mid 1\leq i,j\leq n\}$ subject to the relations making the matrix $(u_{ij})$ orthogonal. The comultiplication is defined by 
\[
\Delta(u_{ij})=\sum_{k=1}^n u_{ik}\tens u_{kj}
\]
and the antipode $S\colon A_0\to A_0$ by $S(u_{ij})=u_{ji}$. These quantum groups were discovered by S. Wang in \cite{wang} and studied further by T.~Banica in \cite{banica-orto}. See also \cite{banica-collins1} and \cite{vaes06}.

\begin{cor}
For $n\geq 3$ we have $\bet_0(A_o(n))=0$.
\end{cor}
\begin{proof}
Denote by $(u_{ij})$ the fundamental corepresentation of $A_o(n)$. Since $S(u_{ij})=u_{ji}$ we have $S^2=\id_{A_0}$ and therefore the Haar state $h$ is tracial (\cite[p.~424]{klimyk}). By \cite[Thm.~7.1]{vaes06} the enveloping von Neumann algebra $\pi_h(A_0)''$ is a $\I\I_1$-factor and  Proposition \ref{faktor-nul} applies.
\end{proof}

\subsection{The commutative case}
	Next we want to investigate the commutative quantum groups. Consider a compact group $G$ and the associated abelian,
compact quantum group $\GG=(C(G),\Delta_c)$. Recall that the comultiplication $\Delta_c\colon C(G)\to C(G)\tens C(G)=C(G\times G)$ is defined by
\[
\Delta_c(f)(s,t)= f(st),
\] 
and that the Haar state and counit are given, respectively, by integration against the Haar probability measure and by evaluation at the identity. In the case when $G$ is 
a connected \emph{abelian} Lie group then $G$ is isomorphic to $\TT^m$ for
some $m\in \NN$ (\cite[Cor.~1.103]{knapp}) and therefore the Pontryagin dual group is $\ZZ^m$. Moreover, the Fourier transform is an isomorphism of quantum groups between $\GG=(C(\TT^m),\Delta_c)$ and
$(C^*_\red(\ZZ^m),\Delta_\red)$. In
particular we have, by Proposition \ref{liggruppe}, that $\bet_0(\GG)=\bet_0(\ZZ^m)=0$ and
\begin{align*}
c_0(\GG) = c_0(\ZZ^m)=\frac{1}{m}=\frac{1}{\dim(G)},
\end{align*}
where the second equality follows from \cite[Thm.~3.7]{LRS}. This motivates the following result.

\begin{thm}\label{dimension-thm}
Let $G$ be a compact Lie group with $\dim(G)\geq 1$ and Haar probability measure $\mu$. Denote by $\GG$ the corresponding
compact quantum group $(C(G),\Delta_c)$. Then $\hh_0(\GG)$ is a finitely
presented and zero-dimensional $L^\infty(G,\mu)$-module, in particular $\bet_0(\GG)=0$, and
\[
c_0(\GG)=\frac{1}{\dim(G)}.
\]
Here $\dim(G)$ is the dimension of $G$ considered as a real manifold.
\end{thm}
For the proof we will need a couple of lemmas/observations probably well known
to most readers. The first one is a purely measure theoretic result.
\begin{lem}\label{ideal-lem}
Let $(X,\mu)$ be measure space and consider $[f_1],\dots,[f_n]\in
L^\infty(X,\RR)$. If we denote by $f$ the function
\[
X\ni x\longmapsto \sqrt{f_1(x)^2 + \dots + f_n(x)^2}\in \RR,
\]
then the ideal $\langle [f_1],\dots,[f_n] \rangle$ in $L^\infty(X,\CC)$ generated by the $[f_i]$'s is equal to the ideal $\langle [f]\rangle$ generated by $[f]$. 
\end{lem}
Here, and in the sequel, $[g]$ denotes the class in $L^\infty(X,\CC)$ of a given function $g$.

\begin{proof}
Consider the real-valued representatives $f_1,\dots,f_n$. Put 
\[
N_i=\{x\in X \ | \ f_i(x)=0\}\qquad \text{ and } \qquad N=\cap_i N_i. 
\]
Note that $N$ is exactly the set of zeros for $f$. \\
\paragraph{''$\mathbf{\subseteq}$''} Let $i\in \{1,\dots, n\}$ be given. We seek $[T]\in L^\infty(X,\CC)$ such that $[f_i]=[T] [f]$. The set $N$ may be disregarded since $f_i$ is zero here. Outside of $N$ we may write
\[
f_i(x)=\frac{f_i(x)}{f(x)} f(x),
\]
and we have $|\frac{f_i(x)}{f(x)}|=\sqrt{\frac{f_i(x)^2}{\sum_j
f_j(x)^2}}\leq 1$. The function
\[
T(x)=
\left\{%
\begin{array}{ll}
   0 & \hbox{if $x\in N$;} \\
    \frac{f_i(x)}{f(x)} & \hbox{if $x\in X\setminus N$,} \\
\end{array}%
\right.
\]
therefore defines a class $[T]$ in $L^\infty(X,\CC)$ with the required properties.\\

\paragraph{''$\supseteq$''} We must find $[T_1],\dots,[T_n]\in L^\infty(X,\CC)$ such that
\begin{align}\label{inkl}
f(x)=T_1(x)f_1(x)+\cdots + T_n(x)f_n(x) \ \textrm{ for
$\mu$-almost all } x\in X.
\end{align}
For any choice of $T_1,\dots, T_n$ both left- and right-hand side
of (\ref{inkl}) is zero when $x\in N$, and it is therefore
sufficient to define $T_1,\dots, T_n$ outside of $N$. Choose a disjoint measurable
partition of $X\setminus N$ into $n$ sets $A_1,\dots, A_n$ such that
\[
|f_k(x)|=\max_i|f_i(x)|>0 \quad \textrm{when } x \in A_k.
\]
Then $1-\chi_{_N}=\sum_{i=1}^n\chi_{_{A_i}}$ and for $x\notin N$ we therefore have
\[
f(x)= \sum_{i=1}^n \chi_{_{A_i}}(x)f(x)=\sum_{i=1}^n
\Big{(}\chi_{_{A_i}}(x)\frac{f(x)}{f_i(x)}\Big{)}f_i(x),
\]
and
\[
|\chi_{_{A_i}}(x)\frac{f(x)}{f_i(x)}|=\chi_{_{A_i}}(x)\sqrt{\frac{\sum_j
f_j(x)^2}{f_i(x)^2}}\leq \sqrt{n}.
\]
Hence the functions $T_1,\dots,T_n$ defined by
\[
T_i(x)=
\left\{%
\begin{array}{ll}
   0 & \hbox{if $x\in N$;} \\
    \chi_{_{A_i}}(x)\frac{f(x)}{f_i(x)} & \hbox{if $x\in X\setminus N$,} \\
\end{array}%
\right.
\]
determines classes $[T_1],\dots, [T_n]$ in $L^\infty(X,\CC)$ with the
required properties.
\end{proof}

\begin{obs}\label{grund-rep}
Every compact Lie group $G$ has a faithful representation  in
$GL_n(\CC)$ for some $n\in \NN$ and for such a
representation $\pi$ it holds that the algebra of all matrix coefficients $C(G)_0$ is generated by the real and imaginary parts of the matrix coefficients of $\pi$. The existence of a faithful representation $\pi$ follows from \cite[Cor.~4.22]{knapp}. Denote by $\pi_{kl}$ its complex matrix coefficients. The fact that $C(G)_0$ is generated by the set 
\[
\{\Re(\pi_{kl}),\Im(\pi_{kl}) \ | \ 1\leq k,l\leq n\}
\]
is the content of \cite[VI, Prop.~3]{chevalley}. 
\end{obs}	 

\begin{obs}\label{kerne-obs}
Let $A$ be a unital $\CC$-algebra generated by elements $x_1,\dots,
x_n$. If $\varps\colon A\to \CC$ is a unital algebra homomorphism then
$\ker(\varps)$ is the two-sided ideal generated by the elements
$x_1-\varps(x_1),\dots,x_n-\varps(x_n)$. This essentially follows from the formula
\[
ab-\varps(ab)=(a-\varps(a))b +\varps(a)(b-\varps(b))
\]
\end{obs}

\begin{obs}\label{exp-obs}

Denote by $\gl_n(\CC)=\mat_n(\CC)$ the Lie
algebra of $GL_n(\CC)$ and by $\exp$ the exponential function
\[
\gl_n(\CC)\ni X\longmapsto \sum_{k=0}^\infty\frac{X^k}{k!}\in
GL_n(\CC),
\]
and consider the map $f\colon \mat_n(\CC)\to \mat_n(\CC)$ given by
$f(X)=\exp(X)-1$. For any norm $\|\cdot\|$ on $\mat_n(\CC)$ there
exist $r,R>0$ and $\lambda_0\in ]0,\frac12]$ such that the following holds: If $X\in \MM_n(\CC)$ has $\|X\|_\infty \leq \frac12$ then for all $\lambda\in [0,\lambda_0]$ we have
\begin{itemize}
\item $\|X\|\leq \lambda \Rightarrow  \|f(X)\|  \leq R\lambda$\\
\item $\|f(X)\|\leq \lambda \Rightarrow \|X\|\leq r\lambda$
\end{itemize}
In the case when the norm in question is the operator norm $\|\cdot\|_\infty$ this is proven, with $\lambda_0=\frac12$ and $R=r=2$, by considering the Taylor expansion around 0 for the scalar versions (i.e.~$n=1$) of $f$ and $f^{-1}$.  Since all norms on finite dimensional vector spaces are equivalent, the general statement follows from this.
\end{obs}
We are now ready to give the proof of Theorem \ref{dimension-thm}.

\begin{proof}[Proof of Theorem \ref{dimension-thm}.]
By Observation \ref{grund-rep}, we may assume that  $G$ is
contained in $GL_n(\CC)$ so that each $g\in G$ can be written as
$g=(x_{kl}(g) + iy_{kl}(g))_{kl}\in GL_n(\CC)$. Again by Observation \ref{grund-rep}, we have that  $A_0\subseteq
A=C(G)$ is given by
\[
A_0=\Alg_\CC(x_{kl},y_{kl} \ | \ 1\leq k,l\leq n),
\]
where $x_{kl}$ and $y_{kl}$ are now considered as functions on $G$. Since $\varps\colon A_0\to \CC$ is given by evaluation at the identity have
\begin{align*}
\varps(x_{kl})& =\varps(y_{kl})=0 \ \textrm{ when } \ k\neq l,\\
\varps(x_{kk})&=\varps(1)=1,\\
\varps(y_{kk})&=0.
\end{align*}
From Observation \ref{kerne-obs} we now get
\[
\ker(\varps ) = \langle x_{kl},y_{kl},x_{kk}-1,y_{kk}  \ | \
1\leq k,l \leq n, k\neq l  \rangle\subseteq A_0.
\]
Thus
\[
\hh_{0}(\GG)=\Tor_0^{A_0}(L^\infty(G),\CC)\simeq
L^\infty(G)\algutens{A_0}\CC \simeq L^\infty(G)/\langle\ker(\varps)\rangle,
\]
where $\langle\ker(\varps)\rangle$ is the ideal in $L^\infty(G)$
generated by $\ker(\varps)\subseteq A_0$. That is, the ideal
\[
\langle x_{kl},y_{kl},1-x_{kk},y_{kk}  \ | \ 1\leq k,l \leq n, k\neq
l  \rangle \subseteq L^\infty(G),
\]
which by Lemma \ref{ideal-lem} is the principal ideal generated by
the (class of the) function
\[
f(g)= \sqrt{\sum_{k,l} (x_{kl}(g)-\delta_{kl})^2 + y_{kl}(g)^2 }.
\]
Note that the zero-set for $f$ consists only of the identity $1\in
G$ and is therefore a null-set with respect to the Haar measure.
Hence we have a short exact sequence
\begin{align}\label{ses}
0\To L^\infty(G)\overset{\cdot f}{\To }L^\infty(G)\To
\hh_{0}(\GG)\To 0.
\end{align}
By additivity of the dimension function (\cite[Thm.~0.6]{Luck98}), this
means that $\bet_0(\GG)=0$. Moreover, the short exact sequence (\ref{ses}) is a finite presentation of $\hh_0(\GG)$ and hence this module has a Novikov-Shubin invariant $\alpha(\hh_0(\GG))$ (in the sense of \cite[Section 3]{Luck97}) which can be computed using the spectral density function
\[
\lambda \longmapsto
h(\chi_{_{[0,\lambda^2]}}(f^2))=\mu(\{g\in G \ | \
f(g)^2\leq \lambda^2\}).
\]
Put $A_\lambda=\{g\in G\mid f(g)^2\leq \lambda^2 \}$. Since the zero-set for $f$ is a $\mu$-null-set we have
\[
\alpha(\hh_0(\GG))=
\left\{%
\begin{array}{ll}
   \displaystyle\liminf_{\lambda\searrow 0}\frac{\ln(\mu(A_\lambda))}{\ln(\lambda)} & \hbox{if $\forall \lambda>0:\mu(A_\lambda)>0$;} \\
    \infty^+ & \hbox{otherwise.} \\
\end{array}%
\right.
\]
Put $m=\dim(G)$ and choose a linear identification of the Lie algebra $\g$ of $G$ with $\RR^m$. By \cite[Thm.~3.31]{warner}, we can choose neighborhoods $V\subseteq \g$ and $U\subseteq G$, around $0$ and $1$ respectively, such that $\exp\colon V\to U$ is a diffeomorphism. This means that 
\[
\varphi=(\exp|_V)^{-1}\colon U\to V\subseteq \g=\RR^m,
\]
constitutes a chart around $1\in G$. Assume without loss of generality that 
\[
V\subseteq \g\cap \{x\in \gl_n(\CC) \ |  \ \|x\|_\infty\leq \tfrac12\}. 
\]
For $g\in G$ we have
\begin{align*}
g\in A_\lambda & \Leftrightarrow \sum_{k,l}(x_{kl}(g)-\delta_{kl})^2 +y_{kl}(g)^2\leq \lambda^2\\
& \Leftrightarrow \|1-g\|_2^2\leq \lambda^2\\
& \Leftrightarrow g \in B_\lambda(1),
\end{align*}
where $B_\lambda(1)$ is the closed $\lambda$-ball in $(\RR^{2n^2},\|\cdot\|_2)$ with center 1. Thus $A_\lambda=G\cap B_\lambda(1)$ and we can therefore choose $\lambda_0\in ]0,\frac12]$ such that $A_{\lambda_0}\subseteq U$. 
Let $\omega$ denote the unique, positive, probability Haar volume
form on $G$ (see e.g.~\cite[Thm.~8.21, 8.23]{knapp}  or \cite[Cor.~15.7]{lee}) and let  $\lambda\in[0,\lambda_0]$. Then
\begin{align*}
\mu(A_\lambda)&= \int_G \chi_{_{A_\lambda}}\d\mu\\
&=\int_{U} \chi_{_{A_\lambda}}\omega\\
&=\int_{V} (\chi_{_{A_\lambda}}\circ \varphi^{-1})(x_1,\dots,x_m) F(x_1,\dots,x_m)\d x_1\cdots \d x_m\\
&= \int_{\varphi(A_\lambda)} F(x_1,\dots,x_m)\d x_1\cdots \d x_m,
\end{align*}
where $F\colon V\to \RR$ is the unique positive function describing
$\omega$ in the local coordinates $(U,\varphi)$. 
By construction we have  $F>0$ on all of $V$ and since $\varphi(A_{\lambda_0})$ is a compact set there exist $C,c>0$ such that
\[
c\leq F(x_1,\dots, x_m)\leq C \ \textrm{ for all } (x_1,\dots, x_m)\in \varphi(A_{\lambda_0}).
\]
For any $\lambda\in[0,\lambda_0]$ we therefore have
\begin{align}\label{lala}
c \nu_m(\varphi(A_\lambda))\leq \mu(A_\lambda)\leq
C\nu_m(\varphi(A_\lambda)),
\end{align}
where $\nu_m$ denotes the Lebesgue measure in $\RR^m=\g$.  
Since $A_\lambda=G\cap B_\lambda(1)$ and 
$\varphi$ is $(\exp|_V)^{-1}$, it follows from Observation \ref{exp-obs} that there exist $d,D>0$ and $\lambda_1\in ]0,\lambda_0]$ such that for all $\lambda\in [0,\lambda_1]$
\[
B_{d\lambda}(0)\cap V \subseteq \varphi(A_\lambda)\subseteq B_{D\lambda}(0)\cap V.
\]
Hence there exist $d',D'>0$ such that for all $\lambda\in [0,\lambda_1]$
\begin{align}\label{dobbeltineq}
d'\lambda^m \leq \nu_m(\varphi(A_\lambda))\leq D' \lambda^m.
\end{align}
From (\ref{lala}) and (\ref{dobbeltineq}) we see that $\mu(A_{\lambda})>0$ for $\lambda\in ]0,\lambda_1]$ and since
\[
\lim_{\lambda\searrow
0}\frac{\ln(d'\lambda^m)}{\ln(\lambda)}=\lim_{\lambda\searrow
0}\frac{\ln(D'\lambda^m)}{\ln(\lambda)}=m,
\]
we also conclude that
\[
\alpha(\hh_0(\GG))=\liminf_{\lambda\searrow 0} \frac{\ln(\mu(A_\lambda))}{\ln(\lambda)}=m=\dim(G).
\]
By definition (\cite[Def.~2.2]{LRS}), the capacity of a finitely presented zero-dimen\-sional module is the inverse of its Novikov-Shubin invariant and thus
\[
c_0(\GG)= c(\hh_0(\GG))=\frac{1}{\dim(G)}.
\]
\end{proof}

\subsection{The finite dimensional case}
In Theorem \ref{dimension-thm} above we only considered compact Lie groups of positive dimension. What is left is the case when $G$ is finite. When $G$ is finite the algebra $C(G)$ is finite dimensional and we have $C(G)_0=C(G)=L^\infty(G)$, which implies vanishing of $H^{(2)}_n(C(G),\Delta_c)$ for $n\geq 1$. For $n=0$ we get
\[
H_0^{(2)}(C(G),\Delta_c)=C(G)\algutens{C(G)}\CC\simeq C(G)\delta_e.
\]
This proves that $H_0^{(2)}(C(G),\Delta_c)$ is a finitely generated projective $C(G)$-module and hence
\[
\bet_0(C(G),\Delta_c)= h(\delta_e)=\int_G \delta_e(g)\d \mu(g)=\frac{1}{|G|}.
\]
Projectivity of $H_0^{(2)}(C(G),\Delta_c)$ implies (cf.~\cite{LRS}) that $c_0(C(G),\Delta_c)=0^-$.\\

This argument generalizes in the following way.

\begin{prop}
Let $\GG=(A,\Delta)$ be a quantum group and assume that $A$ has finite linear dimension $N$. Then
\[
\bet_0(\GG)=\frac{1}{N},
\]
and $\bet_n(\GG)=0$ for all $n\geq 1$. Moreover, $c_n(\GG)=0^-$ for all $n\in \NN_0$.
\end{prop}
\begin{proof}
We first note that for a finite dimensional (hence compact) quantum group the Haar state is automatically tracial (\cite[Thm.~2.2]{vandaele-finite}), so that the numerical $L^2$-invariants make sense. The fact that the higher $L^2$-Betti numbers vanish is a trivial consequence of the fact that $A$ is finite dimensional and therefore equal to both $A_0$ and its enveloping von Neumann algebra. To compute the zeroth $L^2$-Betti number we compute the zeroth $L^2$-homology as
\[
H_0^{(2)}(\GG)= \Tor_0^{A}(A,\CC)\simeq A\algutens{A}\CC = Ae,
\]
where $e$ is the projection in $A$ projecting onto the $ \CC$-summand $A/\ker(\varps)$. Hence
\[
\bet_0(\GG)=\dim_{A}Ae=h(e)=\frac{1}{N},
\]
where the last equality follows, for instance, from \cite[A.2]{woronowicz-pseudo}.\\
Each finite dimensional $C^*$-algebra is a semisimple ring and therefore all modules over it are projective. Hence all capacities of finite dimensional compact quantum groups are $0^-$.
\end{proof}
\vspace{0.3cm}

\section{A vanishing result in the commutative case}\label{comsection}
Throughout this section, $G$ denotes a compact, \emph{connected} Lie group of dimension $m\geq 1$ and $\mu$ denotes the Haar probability measure on $G$. We will also use the following notation:
\begin{align*}
\GG&= (C(G),\Delta_c)\\
A&= C(G)\\
A_0 &= \textrm{The algebra of matrix coefficients}\\
\A&= L^\infty(G,\mu)\\
\U&=\textrm{The algebra of $\mu$-measurable functions on $G$ finite almost everywhere}
\end{align*}
We aim to prove that $\bet_n(\GG)=0$ for all $n\geq 1$. Before doing this, a few comments on the objects defined above. We first note that $\U$ may be identified with the algebra of operators affiliated with $\A$ by \cite[Thm.~5.6.4]{KR1}. In \cite{reich01} it is proved that there is a well defined dimension function $\dim_{\U}(\cdot)$ for modules over $\U$ satisfying properties similar to those enjoyed by $\dim_\M(\cdot)$ (cf.~\cite[Thm.~0.6]{Luck98}). Moreover, by \cite[Thm.~3.1, Prop.~2.1]{reich01} the functor $\U\odot_{\A}-$ is exact and dimension-preserving from the category of $\A$-modules to the category of $\U$-modules. \\
By \cite[Cor.~4.22]{knapp}, we know that $G$ can be faithfully represented in $GL_n(\CC)$ for some $n\in \NN$. Since $GL_n(\CC)$ is a real analytic group (in the sense of \cite{chevalley}), this implies that $G$ has a unique analytic structure making any faithful representation $\pi$ analytic in the following sense: For any $g\in G$ and any function $\varphi$ analytic around $\pi(g)$ the function $\varphi\circ\pi$ is analytic around $g$. This is the content of \cite{chevalley} Chapter IV, \S XIV Proposition 1 and \S XIII Proposition 1. We now choose some fixed faithful representation of $G$ in $GL_n(\CC)$ which will be notationally suppressed in the following. That is, we consider $G$ as an analytic subgroup of $GL_n(\CC)$. Denote by
$\{x_{kl},y_{kl}\mid 1\leq k,l\leq n\}$ the natural $2n^2$ real functions on $GL_n(\CC)$ determining the analytic structure. As noted in Observation \ref{grund-rep}, the algebra $A_0$ is generated by the restriction of the functions $x_{kl}$ and $y_{kl}$ to $G$. Consider some polynomial in the variables $x_{kl}$ and $y_{kl}$; this is clearly an analytic function on $GL_n(\CC)$ and it therefore defines an analytic function on $G$ by restriction. Thus every function in $A_0$ is analytic. \\
The following result is probably well known to experts in Lie theory, but we were unable to find a suitable reference.

\begin{prop}\label{measure-prop}
If $f\in A_0$ is not constantly zero then 
\[
\mu(\{g\in G \ | \ f(g)=0\})=0.
\]
Hence $f$ is invertible in $\U$.
\end{prop}
For the proof we will need the following:
\begin{obs}\label{anal-obs}
Let $V\subseteq \RR^n $ be connected, convex and open and assume that $f\colon V\to \RR$ is analytic. If $f$ is not constantly zero on $V$ then $N=\{x\in V \mid f(x)=0\}$ is a set of Lebesgue measure 0. This is well known in the case $n=1$, since in this case $N$ is at most countable. The general case now follows from this by induction on $n$.
\end{obs}

\begin{proof}[Proof of Proposition \ref{measure-prop}.]
Since $f(x)=0$ iff $\Re(f(x))=\Im(f(x))=0$ we may assume that $f$ is real valued.  Cover $G$ with finitely many precompact, connected, analytic charts 
\[
(U_1,\varphi_1),\dots, (U_t,\varphi_t), 
\]
such that $\varphi(U_i)\subseteq \RR^m$ is convex for each $i\in \{1,\dots, t\}$. Using the local coordinates and the Haar volume form on $G$, it is not hard to see 
that
\begin{align}\label{moved-measure}
\mu(\{g\in U_i\mid f(g)=0\})=0 \ \Leftrightarrow \nu_m(\{x\in \varphi_i(U_i)\mid (f\circ \varphi_i^{-1})(x)=0\})=0.
\end{align}
Here $\nu_m$ denotes the Lebesgue measure in $\RR^m$. Since $f\circ\varphi_i^{-1}$ is analytic it is (by Observation \ref{anal-obs}) sufficient to prove that $f$ is not identically zero on any chart.  
Assume that $f$ is constantly zero on some chart $(U_{i_1},\varphi_{i_1})$. We then aim to show that $f$ is zero on all of $G$, contradicting the assumption. If $G=U_{i_1}$ there is nothing to prove. If not, there exists $i_2\neq i_1$ such that $U_{i_1}\cap U_{i_2}\neq \emptyset$, since otherwise we could split $G$ as the union 
\[
U_{i_1} \cup (\bigcup_{i\neq i_1}U_{i})
\]
of to disjoint, non-empty, open sets, contradicting the fact that $G$ is connected. Since the intersection $U_{i_1}\cap U_{i_2}$ is of positive measure and $f$ is zero on it we conclude, by Observation \ref{anal-obs} and (\ref{moved-measure}), that $f$ is zero on all of $U_{i_2}$ . If $G=U_{i_1}\cup U_{i_2}$ we are done. If not, there exists $i_{3}\notin\{i_1,i_2\}$ such that
\[
U_{i_1}\cap U_{i_3}\neq \emptyset \qquad \textrm{ or } \qquad U_{i_2}\cap U_{i_{3}}\neq \emptyset,
\]
since otherwise $G$ would be the union of two disjoint, non-empty, open sets. In either case we conclude that $f$ is zero on all of $U_{i_3}$. Continuing in this way we conclude that $f$ is zero on all of $G$ since there are only finitely many charts.
\end{proof}
The main result in this section is the following, which should be compared to \cite[Thm.~5.1]{CS}.
\begin{thm}\label{abelian}
Let $Z$ be any $A_0$-module. Then for all $n\geq 1$ we have 
\[
\dim_\A \Tor_n^{A_0}(\A,Z)=0.
\]
\end{thm}
\begin{proof}
As noted in the beginning of this section, we have
\begin{align*}
\dim_{\A} \Tor_n^{A_0}(\A,Z)&= \dim_{\U}(\U\utens{\A} \Tor_n^{A_0}(\A,Z))\\
&=\dim_{\U} \Tor_n^{A_0}(\U,Z).
\end{align*}
We now aim to prove that $\Tor_n^{A_0}(\U,Z)=0$.  For this we first prove the following claim:
\begin{center} 
\emph{Each finitely generated $A_0$-submodule in $\U$ is contained in a finitely generated free $A_0$-submodule.} \\
\end{center}
Let $F$ be a non-trivial, finitely generated $A_0$-submodule in $\U$. We prove the claim by (strong) induction on the minimal number $n$ of generators. If $n=1$ then $F$ is generated by a single element $\varphi\neq 0$, and since all elements in  $A_0\setminus\{0\}$ are invertible in $\U$ (Proposition \ref{measure-prop}) the function $\varphi$ constitutes a basis for $F$. Hence $F$ itself is free. Assume now that the result is true for all submodules that can be generated by $n$ elements, and assume that $F$ is a submodule with minimal number of generators equal to $n+1$. Choose such a minimal system of generators $\varphi_1,\dots,\varphi_{n+1}$. If these are linearly independent over $A_0$ there is nothing to prove. So assume that there exists a non-trivial tuple  $(a_1,\dots, a_{n+1})\in A_0^{n+1}$ such that
\[
a_1 \varphi_1 + \cdots + a_{n+1}\varphi_{n+1}=0,
\]
and assume, without loss of generality, that $a_{1}\neq 0$. Define $F_1$ to be the $A_0$-submodule in $\U$ generated by
\[
a_1^{-1}\varphi_2,\cdots, a_1^{-1}\varphi_{n+1}.
\]
Then $F\subseteq F_1$ and the minimal number of generators for $F_1$ is a most $n$. By the induction hypothesis, there exists a finitely generated free submodule $F_2$ with $F_1\subseteq F_2$ and in particular $F\subseteq F_2$. This proves the claim.\\

	Denote by $(F_i)_{i\in I}$ the system of all finitely generated free $A_0$-submodules in $\U$. By the above claim, this set is directed with respect to inclusion. Since any module is the inductive limit of its finitely generated submodules, the claim also implies that $\U$, as an $A_0$-module, is the inductive limit of the system $(F_i)_{i\in I}$. But, since each $F_i$ is free (in particular flat) and since $\Tor$ commutes with inductive limits we get
\[
\Tor_n^{A_0}(\U, Z)=\lim_{\underset{i}{\to}}\Tor_n^{A_0}(F_i,Z)=0,
\]
for all $n\geq 1$.
\end{proof}
Combining the results of Theorem \ref{abelian} and Theorem \ref{dimension-thm} we get the following.
\begin{cor}\label{abelsk-vaek}
If $G$ is a compact, non-trivial, connected Lie group then 
\[
\bet_n(C(G),\Delta_c)=0, 
\]
for all $n\in \NN_0$.
\end{cor}

\section{Relation to the Connes-Shlyakhtenko approach}\label{cs-section}
In \cite{CS},  A.~Connes and D.~Shlyakhtenko introduced a notion of $L^2$-homology and $L^2$-Betti numbers in the setting of tracial $*$-algebras. More precisely, if $A$ is a weakly dense $*$-subalgebra of a finite von Neumann algebra $\M$ with faithful, normal, trace-state $\tau$, they defined
 (\cite[Def.~2.1]{CS})
\[
H_n^{(2)}(A)=\Tor_n^{A\odot A^\op}(\M\htens\M^\op,A) \quad \textrm{ and } \quad \bet_n(A,\tau)=\dim_{\M\htens\M^\op}H_n^{(2)}(A).
\] 
This generalizes the notion of $L^2$-Betti numbers for groups in the sense that for a discrete group $\Gamma$ we have  $\bet_n(\CC\Gamma,\tau)=\bet_n(\Gamma)$, as proven in \cite[Prop.~2.3]{CS}. In this section we relate the notion of $L^2$-Betti numbers for quantum groups to the Connes-Shlyakhtenko approach. More precisely we prove the following:

\begin{thm}\label{CS-version}
Let $\GG=(A,\Delta)$ be a compact quantum group with tracial Haar state $h$ and algebra of matrix coefficients $A_0$. Then, for all $n\in \NN_0$, we have $\bet_n(\GG)=\bet_n(A_0,h)$, where the latter is the $n$-th $L^2$-Betti number of the tracial $*$-algebra $(A_0,h)$ in the sense of \cite{CS}.
\end{thm}

For the proof of Theorem \ref{CS-version} we will need two small results. Denote by $S\colon A_0\to A_0$ the antipode and by $\varps\colon A_0\to \CC$ the counit. Recall (\cite[p.~424]{klimyk}) that the trace property of $h$ implies that $S^2=\id_{A_0}$ and hence that $S$ is a $*$-anti-isomorphism of $A_0$. Denote by $\M$ the enveloping von Neumann algebra $\pi_h(A_0)''$. In the following we suppress the GNS-representation $\pi_h$ and put $\H=L^2(A,h)$.  Denote by $\bar{\H}$ the conjugate Hilbert space, on which the opposite algebra $A_0^\op$ acts as $a^\op:\bar{\xi}\mapsto \overline{a^*\xi}$.
\begin{lem}\label{implement}
There exists a unitary $V\colon \H\to \bar{\H}$ such that the map 
\[
\B(\H)\supseteq A_0\ni x\overset{\psi}{\longmapsto} (Sx)^\op\in A_0^\op \subseteq \B(\bar{\H})
\]
takes the form $\psi(x)=V x V^*$. In particular, $\psi$ extends to a normal $*$-isomorphism from $\M$ to $\M^\op$.
\end{lem}
 
\begin{proof}
Denote by $\eta$ the inclusion $A_0\subseteq \H$ and note that since $A_0$ is norm dense in $A$ the set $\eta(A_0)$ is dense in $\H$. We now define the map $V$ by
\[
\eta(A_0) \ni \eta(x)\overset{V}{\longmapsto} \overline{\eta(Sx^*)}\in \overline{\eta(A_0)}.
\]
It is easy to see that $V$ is linear and 
\begin{align*}
\|V\eta(x)\|_2^2 &= \|\overline{\eta(Sx^*)}\|_2^2\\
&= \ip{\eta(Sx^*)}{\eta(Sx^*)}\\
&=h((Sx^*)^*S(x^*))\\
&=h(S(x^*x))\\
&=h(x^*x)\\
&=\|\eta(x)\|_2^2,
\end{align*}
and hence $V$ maps the dense subspace $\eta(A_0)$ isometrically onto the dense subspace $\overline{\eta(A_0)}$. Thus, $V$ extends to a unitary which will also be denoted $V$. Clearly the adjoint of $V$ is determined by
\[
\overline{\eta(x)}\overset{V^*}{\longmapsto} \eta(Sx^*).
\] 
To see that $V$ implements $\psi$ we choose some $a\in A_0$ and calculate:
\begin{align*}
\overline{\eta(x)}&\overset{V^*}{\longmapsto} \eta(Sx^*)\\
&\overset{a}{\longmapsto} \eta(aS(x^*))\\
&\overset{V}{\longmapsto}\overline{\eta(S(aS(x^*))^*)}\\
&=\overline{\eta((Sa^*)x)}\\
&=\psi(a)\overline{\eta(x)}.
\end{align*}
\end{proof}
\begin{prop}\label{stjernehomo}
The map $(\id\tens \psi)\circ \Delta\colon A_0\to A_0\odot A_0^\op$ extends to a trace-preserving $*$-homomorphism $\varphi\colon\M\To \M\htens\M^\op$. Here $\psi$ is the map constructed in Lemma \ref{implement} and $\M\htens\M^\op$ is endowed with the natural trace-state $h\tens h^\op$.
\end{prop}
\begin{proof}
The comultiplication is implemented by a multiplicative unitary $W\in \B(\H\htens \H)$ in the sense that
\[
\Delta(a)=W^* (1\tens a)W,
\]
(\cite[page 60]{tuset}) and it therefore extends to a normal $*$-homomorphism, also denoted $\Delta$, from $\M$ to $\M\htens\M$. By Lemma \ref{implement}, the map $\psi\colon \M\to \M^\op$ is normal and therefore $\varphi\colon\M\to \M\htens\M^\op$ is well defined and normal. Since $\varphi$ is normal and $A_0$ is ultra-weakly dense in $\M$ it suffices to see that $\varphi$ is trace-preserving on $A_0$. So, let $a\in A_0$ be given and write $\Delta a=\sum_i x_i\tens y_i\in A_0\odot A_0$. We then have
\begin{align*}
(h\tens h^\op)\varphi(a) &= (h\tens h^\op)(1\tens \psi)(\sum_i x_i\tens y_i)\\
&=(h\tens h^\op)(\sum_i x_i \tens (Sy_i)^\op)\\
&=\sum_i h(x_i)h(y_i)\tag{$h\circ S=h$}\\
&=h(h\tens \id)\Delta(a)\\
&=h(h(a)1_A)\tag{invariance of $h$}\\
&=h(a).
\end{align*}
\end{proof}
We are now ready to give the proof of Theorem \ref{CS-version}.\\

\begin{proof}[Proof of Theorem \ref{CS-version}.]
By Proposition \ref{stjernehomo}, we have that $\varphi$ is a trace-pre\-ser\-ving $*$-homomorphism from $\M$ to $\M\htens\M^\op$. Via $\varphi$ we can therefore consider $\M\htens\M^\op$ as a right $\M$-module and by \cite[Thm.~1.48, 3.18]{sauer-thesis} we have that the functor 
$
(\M\htens\M^\op)\odot_{\M}-
$
is exact and dimension-preserving from the category of $\M$-modules to the category of $\M\htens\M^\op$-modules. Hence 
\begin{align*}
\bet_n(\GG)&= \dim_\M\Tor_n^{A_0}(\M,\CC)\\
&=\dim_{\M\htens\M^\op}(\M\htens\M^\op)\algutens{\M}\Tor_n^{A_0}(\M,\CC)\\
&=\dim_{\M\htens\M^\op}\Tor_n^{A_0}(\M\htens\M^\op, \CC).
\end{align*}
By \cite[Prop.~2.3]{kraehmer} (see also \cite{feng-tsygan}), we have an isomorphism of vector spaces
\begin{align}\label{iso}
\Tor_n^{A_0}(\M\htens\M^\op,\CC)\simeq \Tor_n^{A_0\odot A_0^\op}(\M\htens\M^\op,A_0),
\end{align}

where on the right-hand side $A_0\odot A_0^\op$ acts on $A_0$ in the trivial way and on $\M\htens\M^\op$ via the natural inclusion $\M\htens\M^\op \supseteq A_0\odot A_0^\op$. This isomorphism respects the natural left action of $\M\htens\M^\op$, since on both sides of (\ref{iso}) only the multiplication from the right on $\M\htens\M^\op$ is used to compute the $\Tor$-groups. The right-hand side of (\ref{iso}) is, by definition, equal to the $L^2$-homology of $A_0$ in the sense of \cite{CS} and the statement follows.
\end{proof}

\begin{cor}
Let $G$ be a non-trivial, compact, connected Lie group with Haar measure $\mu$ and denote by $A_0$ the algebra of matrix coefficients arising from irreducible representations of $G$. Then, for all $n\in \NN_0$, we have $\bet_n(A_0, \d \mu)=0$.
\end{cor}
\begin{proof}
This follows from Theorem \ref{CS-version} and Corollary \ref{abelsk-vaek} in conjunction
\end{proof}


\end{document}